\documentclass[11pt,a4paper]{article}
\usepackage[utf8]{inputenc}
\usepackage{amsmath}
\usepackage{amsfonts}
\usepackage{amssymb}
\usepackage{makeidx}
\usepackage{graphicx}
\usepackage{mathabx} 
\usepackage{dsfont} 
\usepackage[left=2cm,right=2cm,top=2cm,bottom=2cm]{geometry}
\usepackage{color}
\usepackage{verbatim} 
\usepackage{tocbibind} 
\usepackage{hyperref} 
\usepackage{url}
\usepackage{mathtools} 
\usepackage{dsfont}
\usepackage{ulem}

\newtheorem{proposition}{Proposition}[section]
\newtheorem{definition}{Definition}[section]
\newtheorem{theorem}{Theorem}[section]
\newtheorem{remark}{Remark}[section]

\begin{document}

\setlength{\parindent}{0pt} 
\begin{center}
\Large\textbf{THE COMPRESSIBLE NAVIER-STOKES EQUATIONS WITH SLIP BOUNDARY CONDITIONS OF FRICTION TYPE} \\[10mm]
\large{\v{S}\' ARKA NE\v{C}ASOV\' A}$^1$, \large{JUSTYNA OGORZALY}$^{1,3}$, \large{JAN SCHERZ}$^{1,2}$ \\[10mm]
\end{center}

\begin{itemize}
\item[$^1$] Institute of Mathematics of the Czech Academy of Sciences, \v{Z}itn\' a 25, Prague 1, 11567
\item[$^2$] Department of Mathematical Analysis, Faculty of Mathematics and Physics, Charles University in Prague, Sokolovsk\' a 83, Prague 8, 18675
\item[$^3$] Institute of Mathematics and Cryptology, Military University of Technology, ul. Gen. Sylwestra Kaliskiego 2,
00-908 Warsaw 
\end{itemize}
{\bf{Abstract.}} {We study a mathematical model of a viscous compressible fluid obeying the slip boundary condition of friction type. We present a notion of weak solutions to this model, in which the momentum equation and the associated energy inequality are combined into a single relation. Moreover, the slip boundary condition of friction type is incorporated into this relation by the use of a boundary integral. Our main result proves the existence of such weak solutions. The proof of this result combines the classical existence theory for the compressible Navier-Stokes equations with an approximation of the aforementioned boundary integral via a convex regularization of the absolute value function.} 
\section{Introduction}
The subject of this article is a new model of a compressible viscous fluid satisfying the so-called \textit{slip boundary condition of friction type (SBCF)}. In this paper we prove the existence of weak solutions to this model.


{\it Classical} (or strong) solutions to the compressible Navier-Stokes equations can be expected only for small data (or, more generally, for data close to an equilibrium). The first such result for the Cauchy problem for the Navier-Stokes-Fourier system (when heat conductivity is included) goes back to the eighties (\cite{MaNi}; for further developments, see e.g. \cite{MuZa} or \cite{Da}. However, classical solutions are not known to exist globally in time if the data is arbitrary. The concept of {\it weak solutions} was for the first time successfully used by Lions (see \cite{Lions2}) in the case of isentropic flow. In this book, several kinds of boundary conditions were considered: The \textit{no-slip boundary condition}, which describes the vanishing of the fluid velocity on the boundary of the domain, periodic boundary conditions as well as the case of a fluid covering the whole space. A detailed proof of the existence of weak solutions in the case of the no-slip boundary condition can further be found in \cite{novotnystraskraba}. A weak solution in the case of heat-conducting fluids satisfying the no-slip boundary condition was for the first time constructed by Feireisl (see \cite{feireisldynamics}) by combining the internal energy balance and the global energy balance. Another approach, presented by Feireisl and Novotn\'y, is based on the entropy inequality (see \cite{singularlimits}). In the latter book the case of the {\it complete-slip boundary condition}, i.e.\ the case of fluids for which the normal component of the velocity vanishes on the boundary, is additionally taken into consideration. Moreover, the existence of weak solutions in the case of incompressible fluids is treated for example in \cite{lions} for the no-slip boundary condition, periodic boundary conditions as well as in the whole space $\mathbb{R}^N$.

The no-slip boundary condition has been the most widely used given its success in reproducing the standard velocity profiles for incompressible/compressible viscous fluids for many years. The no-slip hypothesis seems to be in good agreement with experiments but it can lead to certain rather surprising conclusions e.g. the most striking one being the absence of collisions of rigid objects immersed in a linearly viscous fluid \cite{HES,HIL}.

The Navier-Stokes equations have also been studied in combination with more uncommon boundary conditions. The so-called {\it Navier boundary condition}, which allows for slip, offers more freedom and is likely to provide a physically acceptable solution at least to some of the paradoxical phenomena resulting from the no-slip boundary condition, see, e.g.\ Moffat \cite{MOF}. Recent developments in macrofluidic and nanofluidic technologies have renewed interest in the slip behavior that may become significant in the small spatial scales even for a relatively small Reynolds number (cf. Priezjev and Troian \cite{PRTR}). Mathematically, the behavior of the tangential component of the velocity is a delicate issue.

We further mention the Coulomb friction law boundary condition, which is used for the description of fluids that can slip on the boundary provided that the tangential component of the stress tensor is sufficiently large. In \cite{bmt} the existence of weak solutions to the incompressible Navier-Stokes equations satisyfing this boundary condition is proved in the case of two and three spatial dimensions. Another boundary condition modelling this phenomenon is the slip boundary condition of friction type introduced by H. Fujita in \cite{fuj1} and \cite{fuj2} H. Fujita for the stationary Stokes and Navier-Stokes equations. The same boundary condition was studied for the incompressible Navier-Stokes equations in \cite{kas2}, wherein the existence of solutions is proved globally in time in the $2$D case and locally in time in the $3$D case. A numerical analysis of the slip boundary condition of friction type can be found e.g.\ in \cite{kas1}. Moreover, some applications to real world problems with numerical simulations are given in \cite{kas3}, \cite{HM} and \cite{kas4}.

In the present article we combine, for the first time, the compressible Navier-Stokes equations with the slip boundary condition of friction type. We prove the existence of weak solutions in this setting. Since the slip boundary condition of friction type is particularly interesting for the modelling of fluids in moving domains or fluid-structure interaction, our result can be considered as a first stepping stone towards the study of these more sophisticated problems. From the mathematical point of view, the main novelty in our existence proof lies in the addition of one further approximation level to the classical approximation method used for the construction of weak solutions to the compressible Navier Stokes equations with the no-slip boundary condition, c.f.\ for example \cite{novotnystraskraba}. This additional approximation level follows closely the approximation methods used in \cite{bmt} and \cite{kas2} in the case of the incompressible Navier-Stokes equations with the Coulomb friction law boundary condition and the slip boundary condition of friction type respectively. It consists of the addition of a boundary integral to the momentum equation which contains the gradient of a smooth and convex approximation of the absolute value of the velocity field, c.f.\ \eqref{weakmomeq} below. Due to the convexity of the approximation, this boundary integral can later be replaced by the desired boundary integral which expresses the slip boundary condition of friction type in our weak formulation \eqref{momentumequation}. Another novelty results from the fact that, as in the incompressible case in \cite{kas2}, the weak formulation to our problem merges the momentum equation and the energy inequality into one single relation, c.f.\ \eqref{momentumequation} below. For technical reasons, however, we also need to study the momentum equation separately in order to deduce the same improved density estimates and the effective viscous flux identity as in the existence proof in the case of the no-slip boundary condition in \cite{novotnystraskraba}, which are required for passing to the limit in the pressure term. As a consequence we are forced to also pass to the limit in the momentum equation separately on every approximation level, which leads to a relation which we refer to as the alternative momentum equation, c.f.\ Remark \ref{remarkaltmomeq} below.

The paper is organised as follows. In Section \ref{model} we present the full model. A corresponding weak formulation of this model is presented in Section \ref{weaksolutions}. In the same section we further show that this weak formulation constitutes a suitable definition of weak solutions and present our main result. The full proof of the main result extends across the Sections \ref{solapproxprob}--\ref{alphalim}.

\section{Model}\label{model}
The model which we study in this paper is as follows. We consider a viscous compressible fluid occupying an open and bounded domain $\Omega \subset \mathbb{R}^3$ with locally Lipschitz boundary $ \Gamma = \partial \Omega$ and with outward unit normal vector $\text{n}$ on $\Gamma$. The density $\rho \colon (0,T) \times \Omega \to \mathbb{R}$ and the velocity field $u \colon (0,T) \times \Omega \to \mathbb{R}^3$ of the fluid are determined via the system
\begin{align}
\partial _t \rho + \nabla \cdot (\rho u) = 0 \quad \quad &\text{in } (0,T)\times \Omega, \label{1} \\
\partial _{t} (\rho u) + \nabla \cdot (\rho u \otimes u) = \nabla \cdot \sigma + \rho f \quad \quad &\text{in } (0,T)\times \Omega, \label{2} \\
\rho(0) = \rho^0,\quad (\rho u)(0) = q \quad \quad &\text{in } \Omega, \label{3} \\
u\cdot \text{n} = 0 \quad \quad &\text{on } (0,T)\times \Gamma, \label{4} \\
\left| \left(\sigma \text{n} \right)_\tau \right| \leq g,\quad \left( \sigma \text{n} \right)_\tau \cdot u_\tau + g \left| u_\tau \right| = 0 \quad \quad &\text{on } (0,T)\times \Gamma, \label{5}
\end{align}
where the Cauchy stress tensor
\begin{align} \label{cs}
\sigma = \sigma(u,p) := 2\nu \mathbb{D}(u) + \lambda (\nabla \cdot u)\text{Id}  - p \text{Id},\quad \quad \mathbb{D}(u) := \frac{1}{2} \nabla u + \frac{1}{2} (\nabla u)^T \nonumber
\end{align}

with the viscosity coefficients $\nu,\lambda \in \mathbb{R}$, satisfying 
\begin{align}
\nu > 0,\quad \quad \lambda + \nu \geq 0, \nonumber
\end{align}

can be split into its normal component $\sigma_\text{n} = \sigma \text{n} \cdot \text{n}$ and its tangential component $\sigma_\tau= \sigma \text{n} - \sigma_\text{n}\text{n}$. Moreover the positive constant $g$ in \eqref{5} is the threshold of slippage. Further the pressure $p$ is defined by the isentropic constitutive relation
\begin{equation} \label{pres}
p = a\rho^\gamma,\quad \gamma > \frac{3}{2},\ a > 0. \nonumber
\end{equation}
In the considered model the equations \eqref{1} and \eqref{2} represent the continuity equation and the momentum equation, respectively. The initial conditions are presented in \eqref{3}. Finally the equations \eqref{4} and \eqref{5} represent the slip boundary condition of friction type.
\section{Weak formulation and main result} \label{weaksolutions}
Here we present the definition of a weak solution to the system \eqref{1}--\eqref{5} and state our main result. To this end we denote by $H_\text{n}^1(\Omega; \mathbb{R}^3)$ the Sobolev space of all functions in $H^1(\Omega;\mathbb{R}^3)$ whose normal component vanishes on the boundary,
\begin{align}
H_\text{n}^1\left(\Omega;\mathbb{R}^3 \right) := \left\lbrace u \in H^1\left(\Omega;\mathbb{R}^3 \right):\ u\cdot \text{n} = 0 \ \text{on } \Gamma \right\rbrace. \nonumber
\end{align}

\begin{definition} \label{weaksoldef}
Let $T>0$ and let $\Omega \subset \mathbb{R}^3$ be a bounded domain. Let $\nu, \lambda, a, \gamma \in \mathbb{R}$ be given constants which satisfy
\begin{align}
\nu, a > 0,\quad \gamma > \frac{3}{2},\quad \nu + \lambda \geq 0. \label{assumptionsconstants}
\end{align}

Further assume that $f \in L^\infty ((0,T)\times \Omega)$, $g \in L^2((0,T)\times \Gamma)$ and assume the initial data to satisfy the conditions
\begin{align}
0 \leq \rho_0 \in L^\gamma (\Omega),\quad q \in L^1(\Omega),\quad \frac{\left|q\right|^2}{\rho_0} \in L^1(\Omega),\quad q = 0 \ \ \text{a.e. in } \left\lbrace x \in \Omega:\ \rho_0(x) = 0 \right\rbrace \label{initialdata}
\end{align}
 Then a pair of functions $(\rho, u)$, such that
\begin{align}
0 \leq \rho \in L^\infty \left(0,T;L^\gamma \left(\Omega;\mathbb{R}\right) \right) \bigcap C\left([0,T];L^1\left(\Omega;\mathbb{R} \right) \right) \quad \mbox{and} \quad u \in L^2\left(0,T;H_{\operatorname{n}}^1\left(\Omega, \mathbb{R}^3 \right) \right), \label{regularityweaksol}
\end{align}

is said to be a weak solution to the system \eqref{1}--\eqref{5} if it satisfies:
\begin{itemize}
\item[(i)] the continuity equation in the distributional sense,
\begin{align}
\partial_t \rho + \nabla \cdot \left( \rho u \right) = 0 \quad &\text{in } \mathcal{D}'\left((0,T) \times \mathbb{R}^3 \right), \label{-214}
\end{align}

and the renormalized sense,
\begin{align}
&\ \partial_t \zeta (\rho) + \nabla \cdot \left( \zeta \left(\rho \right)u \right) + \left[ \zeta'\left(\rho \right)\rho - \zeta \left(\rho \right) \right] \nabla \cdot u = 0 \quad \text{in } \mathcal{D}' \left((0,T) \times \mathbb{R}^3 \right), \label{-215} \\
\text{for all} \quad \zeta \in &C^1\left( [0,\infty) \right): \quad \left|\zeta'(r)\right| \leq cr^{\sigma} \quad \forall r \geq 1 \quad \text{for certain } c>0,\ \sigma > -1, \label{-216}
\end{align}

\item[(ii)]the momentum and energy inequality
\begin{align}
&\int_\Omega \frac{1}{2}\rho(0) \left|u(0)\right|^2 + \frac{a\rho^\gamma(0)}{\gamma - 1} \ dx - \int_\Omega \frac{1}{2} \rho(\tau) \left| u(\tau) \right|^2 + \frac{a\rho^\gamma (\tau)}{\gamma - 1} \ dx \nonumber \\
+& \int_0^\tau \int_\Omega - \rho u \cdot \partial_t \phi - (\rho u \otimes u) : \nabla \phi + \left(2\nu \mathbb{D}(u) + \lambda \left( \nabla \cdot u \right)\operatorname{id} \right) : \left[\nabla \phi - \nabla u \right] \nonumber \\
&-p \operatorname{id} : \nabla \phi - \rho f \cdot \left[ \phi - u \right] \ dxdt + \int_0^\tau \int_{\partial \Omega} g |\phi| - g |u| d\Gamma dt \geq 0 
\label{momentumequation}
\end{align}
for almost all $\tau \in [0,T]$ and all $\phi \in \mathcal{D}((0,\tau)\times \overline{\Omega})$ with $\phi \cdot \text{n}|_{\partial \Omega} = 0$ and
\item[(iii)] the initial conditions
\begin{align}
\rho(0) = \rho_0,\quad \quad \lim_{\tau \rightarrow 0+} \int_\Omega \rho(\tau,x)u(\tau,x) \cdot \phi(x) \ dx = \int_\Omega q(x) \cdot \phi (x)\ dx \label{initialcond}
\end{align}

for all $\phi \in \mathcal{D}(\overline{\Omega})$ with $\phi \cdot \operatorname{n}|_{\partial \Omega} = 0$.
\end{itemize}

\end{definition}

In order to make sure that Definition \eqref{weaksoldef} is a suitable definition of weak solutions, we show that any classical solution to the system \eqref{1}--\eqref{5} is also a weak solution and, vice versa, any weak solution with a sufficient amount of regularity solves the problem \eqref{1}--\eqref{5} in the classical sense. In order to obtain the variational inequality \eqref{momentumequation} from the system \eqref{1}--\eqref{5} we first pick an arbitrary time $\tau \in [0,T]$ and multiply the momentum equation \eqref{2} by an arbitrary function 
$\phi \in \mathcal{D}((0,\tau) \times \overline{\Omega})$ with $\phi \cdot \operatorname{n}|_{\partial \Omega} = 0$. Integrating (by parts) over $(0,\tau)\times \Omega$ we obtain the identity
\begin{align}
\int_0^\tau \int_\Omega - \rho u \cdot \partial_t \phi - (\rho u \otimes u) : \nabla \phi\ + \sigma : \nabla \phi \ dxdt =\int_0^\tau \int_\Omega \rho f \cdot \phi \ dxdt + \int_0^\tau \int_{\partial \Omega} (\sigma \text{n})_\tau \cdot \phi\ d\Gamma dt. \label{momentumasusual}
\end{align}

Similarly, we test the momentum equation \eqref{2} by $u$ and subtract from it the continuity equation \eqref{1} tested by $\frac{1}{2}|u|^2$. Hence we infer the energy inequality
\begin{align}
&\frac{1}{2} \int_\Omega \rho (\tau) |u(\tau)|^2 + \frac{a\rho^\gamma (\tau)}{\gamma - 1} \ dx + \int_0^\tau 2\nu \left| \mathbb{D}(u) \right|^2 + \lambda \left| \nabla \cdot u \right|^2\ dxdt \nonumber \\
=& \frac{1}{2} \int_\Omega \rho(0) \left| u(0) \right|^2 + \frac{a\rho^\gamma(0)}{\gamma - 1} \ dx + \int_{0}^\tau \int_\Omega \rho f \cdot u\ dxdt + \int_0^\tau \int_{\partial \Omega} \left( \sigma \text{n} \right)_{\tau} \cdot u \ d\Gamma dt. \nonumber
\end{align}

The last identity we substract from the equation \eqref{momentumasusual} and then we use the boundary condition \eqref{5}. Finally we obtain
\begin{align}
&\int_\Omega \frac{1}{2}\rho(0) \left|u(0)\right|^2 + \frac{a\rho^\gamma(0)}{\gamma - 1} \ dx - \int_\Omega \frac{1}{2} \rho(\tau) \left| u(\tau) \right|^2 + \frac{a\rho^\gamma (\tau)}{\gamma - 1} \ dx + \int_0^\tau \int_\Omega - \rho u \cdot \partial_t \phi - (\rho u \otimes u) : \nabla \phi\  \nonumber \\
&\ + \left(2\nu \mathbb{D}(u) + \lambda \left( \nabla \cdot u \right)\operatorname{id} \right) : \left[\nabla \phi - \nabla u \right] -p \operatorname{id} : \nabla \phi - \rho f \cdot \left[ \phi - u \right] \ dxdt + \int_0^\tau \int_{\partial \Omega} g |\phi| - g |u| d\Gamma dt \nonumber \\
=& \int_0^\tau\int_{\partial \Omega} (\sigma \text{n})_\tau \cdot \left[\phi - u \right]\ d\Gamma dt + \int_0^\tau \int_{\partial \Omega} g |\phi| - g |u| d\Gamma dt \label{equality} \\
=& \int_0^\tau\int_{\partial \Omega} (\sigma \text{n})_\tau \cdot \phi \ d\Gamma dt + \int_0^\tau \int_{\partial \Omega} g |\phi| d\Gamma dt \geq 0, \label{inequality}
\end{align} 

which is exactly the variational inequality \eqref{momentumequation}. Conversely, we need to check that any sufficiently regular weak solution in the sense of Definition \ref{weaksolutions} also satisfies the system \eqref{1}--\eqref{5} in the classical sense. Hereof, the continuity equation \eqref{1}, the initial condition \eqref{3} and the boundary condition \eqref{4} are clear. For the derivation of the momentum equation we test the variational inequality \eqref{momentumequation} by $\psi u \pm \phi$ for some arbitrary functions $\psi \in \mathcal{D}(0,T)$, $\phi \in \mathcal{D}((0,T)\times \Omega)$. Under exploitation of the assumed smoothness of $\rho$ and $u$ this yields the relation
\begin{align}
&\int_\Omega \frac{1}{2}\rho(0) \left|u(0)\right|^2 + \frac{a\rho^\gamma(0)}{\gamma - 1} \ dx - \int_\Omega \frac{1}{2} \rho(\tau) \left| u(\tau) \right|^2 + \frac{a\rho^\gamma (\tau)}{\gamma - 1} \ dx + \int_0^T \int_\Omega \partial_t \left( \rho u \right) \cdot \left(\psi u \pm \phi\right) \nonumber \\
&\ + \left( \nabla \cdot \left( \rho u \otimes u \right)\right) \cdot \left(\psi u \pm \phi \right) + \left( \nabla \cdot \left(2\nu \mathbb{D}(u) + \lambda \left( \nabla \cdot u \right)\operatorname{id} \right)\right) \cdot \left[ \left( 1 - \psi \right)u \pm \phi \right] + \nabla p \cdot \left( \psi u \pm \phi \right) \nonumber \\
&\ + \rho f \cdot \left[ \left(1 - \psi \right) u \pm \phi \right] \ dxdt + \int_0^T \int_{\partial \Omega} g |\psi u| - g |u| d\Gamma dt \geq 0 \nonumber
\end{align}

Letting $\psi \rightarrow 1$ we see that the $\phi$-independent terms in this inequality cancel each other and we are left with the equality
\begin{align}
&\int_0^T \int_\Omega \partial_t \left(\rho u \right) \cdot \phi + \left( \nabla \cdot (\rho u \otimes u)\right) \cdot \phi\ - \left( \nabla \cdot \left(2\nu \mathbb{D}(u) + \lambda \left( \nabla \cdot u \right)\text{id} \right)\right) \cdot \phi + \nabla p \cdot \phi - \rho f \cdot \phi \ dxdt = 0. \label{49}
\end{align}

Hence, by the arbitrary choice of $\phi \in \mathcal{D}((0,T)\times \Omega)$, the classical formulation \eqref{2} of the momentum equation is satisfied. This in particular implies that the identity \eqref{equality} again holds true. Subtracting \eqref{equality} from the given variational inequality \eqref{momentumequation} we find the estimate
\begin{align}
-\int_0^\tau\int_{\partial \Omega} (\sigma \text{n})_\tau \cdot \left[\phi - u \right]\ d\Gamma dt \leq \int_0^\tau \int_{\partial \Omega} g |\phi| - g |u| d\Gamma dt \label{48}
\end{align}

for any $\phi \in \mathcal{D}((0,\tau)\times \overline{\Omega})$ with $\phi \cdot \text{n}|_{\partial \Omega} = 0$ and, by a density argument, for any $\phi \in L^2(0,\tau;H_{\text{n}}^1(\Omega))$. We choose $s=T$ and test this inequality by $u \pm \phi$. Hence, replacing $\phi$ by $u \pm \phi$, we conclude, from the reverse triangle inequality, the estimate
\begin{align}
\left| \int_0^\tau \int_{\partial \Omega} (\sigma \text{n})_\tau \cdot \phi \ d\Gamma dt \right| \leq \int_0^\tau \int_{\partial \Omega} g |\phi| d\Gamma dt. \nonumber
\end{align}

It follows that $|(\sigma \text{n})_\tau|\leq g$ on $\partial \Omega$ and further, from the estimate \eqref{48} with the choice $\phi = 0$, that $(\sigma \text{n})_\tau \cdot u + g |u| = 0$ on $\partial \Omega$. Hence also the boundary condition \eqref{5} is satisfied. \\

We are now in the position to present the main result of our article, which is as follows:

\begin{theorem}
\label{mainresult}
Let $T > 0$ and let $\Omega \subset \mathbb{R}^3$ be a bounded domain of class $C^{2,\eta} \bigcup C^{0,1}$ for some $\eta > 0$. Let the data $\nu, \lambda, a, \gamma \in \mathbb{R}$, $f \in L^\infty((0,T)\times \Omega)$, $g \in L^2((0,T)\times \partial \Omega)$, $\rho_0 \in L^\gamma (\Omega)$ and $q \in L^1(\Omega)$ satisfy the conditions \eqref{assumptionsconstants}--\eqref{initialdata}. Then there exists a weak solution $(\rho, u)$, in the sense of Definition \ref{weaksoldef}, to the system \eqref{1}--\eqref{5}.
\end{theorem}

We remark that the $C^{2,\eta}$-regularity of $\Omega$ in Theorem \ref{mainresult} is necessary for the construction of the density in the approximate system in Section \ref{solapproxprob} below, c.f.\ \cite[Lemma 3.1, Theorem 10.22, Theorem 10.23]{singularlimits}, \cite[Proposition 7.39]{novotnystraskraba}. Moreover, the $C^1$-regularity of $\Omega$ is needed to extend $u$ to an $L^2(0,T;H^1(\mathbb{R}^3))$-function when showing that the couple $(\rho, u)$ satisfies the renormalized continuity equation in Section \ref{epslimit} below, c.f.\ \cite[Section 5.4, Theorem 1]{evans}.

\section{Approximate system}\label{aprox}
In this section we present an approximate version of the problem introduced in Section \ref{model}, followed by a brief explanation of the individual approximation levels. We fix four parameters $n \in \mathbb{N}$, $\delta, \epsilon, \alpha >0$, each of them associated to one of these approximation levels. We further fix some parameter $\beta > \max \lbrace \gamma, 4 \rbrace$. By $V_n \subset C^2(\overline{\Omega};\mathbb{R}^3) \subset L^2(\Omega;\mathbb{R}^3)$ we denote an $n$-dimensional vector space equipped with the $L^2(\Omega)$-inner product, such that
\begin{align}
\bigcup _{n\in \mathbb{N}}V_n \ \ \text{is dense in} \ \ W_{\text{n}}^{1,p}(\Omega) := \left\lbrace \phi \in W^{1,p}(\Omega):\ \left. \phi \cdot \text{n} \right|_{\Gamma} = 0 \right\rbrace \quad \quad \forall 1 \leq p < \infty. \label{galerkindensity}
\end{align}

For technical reasons (c.f. the deduction of the convergence \eqref{cweakmom} below) we assume that without loss of generality the sequence of spaces $(V_n)_n$ contains a subsequence $(V_{0,n})_n$ of spaces $V_{0,n}\subset C_0^2(\overline{\Omega};\mathbb{R}^3)$ such that
\begin{align}
\bigcup _{n\in \mathbb{N}}V_{0,n} \ \ \text{is dense in} \ \ W_0^{1,p}(\Omega) \quad \quad \forall 1 \leq p < \infty. \label{galerkindensity2}
\end{align}

Moreover, following the approximation method used for the proof of the existence of weak solutions to the Navier-Stokes equations with the Coulomb friction law boundary condition in \cite[Section 3]{bmt}, we denote by
\begin{align}
j_\delta (v) := \left\{
\begin{matrix}
|v| & \text{for } |v| > \delta, \\
\frac{|v|^2}{2\delta} + \frac{\delta}{2} & \text{for } |v| \leq \delta,
\end{matrix}\right. \nonumber
\end{align}

a convex approximation $j_\delta \in C^1(\mathbb{R}^3)\bigcap C^{1,1}_{\text{loc}}(\mathbb{R}^3)$ of the absolute value function. We remark that while in \cite{bmt} the local Lipschitz-continuity of the gradient of $j_\delta$ is not required, it is necessary in our setting in order to achieve continuity of the operator $\mathbb{T}$ in the fixed point argument for the construction of an approximate solution in Section \ref{solapproxprob} below. The approximation $j_\delta$ further has the properties
\begin{align}
j_\delta (0) =& 0, \label{jorigin}\\
\operatorname{grad} j_\delta(v) \cdot v \geq & 0 \quad \forall \ v \in \mathbb{R}^3, \label{nonnegj}\\
\left| \operatorname{grad} j_\delta(v) \right| \leq & 1 \quad \forall \ {v} \in \mathbb{R}^3, \label{conditionj}\\
\left| j_\delta(v) - \left| v \right| \right| \leq & \delta \quad \forall \ v \in \mathbb{R}^3, \label{convergencej}
\end{align}

where $\operatorname{grad} j_\delta$ denotes the gradient of $j_\delta$. Our approximate problem on the highest approximation level consists of finding functions
\begin{align}
\rho_\delta \in& W:= \left\lbrace \psi \in C\left([0,T]; C^{2,\eta} \left( \overline{\Omega} \right) \right) \bigcap C^1\left([0,T]; C^{0,\eta} \left( \overline{\Omega} \right) \right):\ \left. \nabla \psi \cdot \text{n} \right|_{\Gamma} = 0 \right\rbrace, \label{approxrho} \\
u_\delta \in& C\left([0,T];V_n \right) \label{approxu}
\end{align}

which satisfy the approximate continuity equation
\begin{align}
\partial_t \rho_\delta + \nabla \cdot \left( \rho_\delta u_\delta \right) =& \epsilon \Delta \rho_\delta, \quad \quad \left. \nabla \rho_\delta \cdot \text{n} \right|_{\Gamma} = 0. \label{weakconteq}
\end{align}

in $(0,T)\times \Omega$ and the approximate momentum equation
\begin{align}\label{weakmomeq}
\int_\Omega \partial_t \left(\rho_\delta u_{\delta} \right) \cdot \phi \ dx =& \int _\Omega \left( \rho_\delta u_\delta \otimes u_{\delta} \right): \nabla \phi - 2\nu \mathbb{D}\left(u_{\delta}\right):\mathbb{D} (\phi) - \lambda (\nabla \cdot u_\delta ) (\nabla \cdot \phi) + a\rho_\delta ^\gamma \nabla \cdot \phi + \alpha \rho_\delta ^\beta \nabla \cdot \phi \nonumber \\
&+ \rho_\delta f\cdot \phi - \epsilon \left( \nabla u_{\delta} \nabla \rho_\delta \right) \cdot \phi\ dx - \int_{\partial \Omega} g \operatorname{grad} j_\delta(u_\delta) \cdot \phi\ d\Gamma
\end{align}

in $[0,T]$ for all $\phi \in C([0,T];V_n)$ as well as the initial conditions
\begin{align}
\rho(0,x) = \rho_0(x),\quad u(0,x) = u_0(x) \quad \quad \forall x \in \Omega. \label{initialcondapproxi}
\end{align}

Here the initial data $u_0$ for the velocity field is defined by
\begin{align}
u_0 := P_n \left( \frac{q}{\rho_0} \right) \in V_n, \label{indata}
\end{align}

where $P_n$ denotes the orthogonal projection from $L^2(\Omega)$ onto $V_n$, and the initial data $\rho_0$, $q$ is assumed to satisfy the additional regularity criteria
\begin{align}
\rho _0 \in C^{2,\eta}\left(\overline{\Omega}\right),\quad \quad 0 < \alpha \leq \rho_0 \leq \alpha^{-\frac{1}{2\beta}},\quad \quad \left. \nabla \rho_0 \cdot \operatorname{n} \right|_{\Gamma} = 0,\quad \quad q \in C^2\left(\overline{\Omega}\right). \label{-342}
\end{align}

Having introduced the full approximate problem \eqref{approxrho}--\eqref{initialcondapproxi} we now give a short explanation of the individual approximation levels in the order, in which we will later pass to the limit in them, beginning with the $\delta$-level. On this level, following the proof of the existence of weak solutions to the incompressible Navier-Stokes equations with the Coulomb friction law boundary conditions in \cite{bmt}, we add a boundary integral containing the quantity $\operatorname{grad} j_\delta (u_\delta)$ to the momentum equation. The convexity of $j_\delta$ then allows us to transform the approximate momentum equation \eqref{weakmomeq} into an inequality in which the desired boundary condition is incorporated in the same way as in the momentum and energy inequality \eqref{momentumequation} in our weak formulation. \\
The remaining approximation levels coincide precisely with the corresponding approximation levels in the classical theory of the existence of weak solutions to the compressible Navier-Stokes equations, which can be found for example in \cite[Chapter 7]{novotnystraskraba}. On the $n$-level we carry out a Galerkin approximation, which allows us to find a solution to the approximate system. More precisely, this procedure reduces the problem to a finite dimensional problem in the spatial component, which can be solved via the classical theory of ordinary differential equations and a fixed point argument. The reason why we pass to the limit with respect to $\delta \rightarrow 0$ before passing to the limit in the Galerkin approximation lies in the high spatial regulariy available on the Galerkin level. This regularity allows us to achieve uniform convergence of the velocity field $u_\delta$ when letting $\delta$ tend to zero, which is required for passing to the limit in the quantity $j_\delta (u_\delta)$. \\
On the $\epsilon$-level the additional quantity $\epsilon \Delta \rho_\delta$ is added to the continuity equation. This procedure (c.f.\ \cite[Section 7.6]{novotnystraskraba}), known as the parabolic regularization of the continuity equation, is required to make sure that the density in our approximate system and consequently also in our final system is non-negative. For the sake of preserving an energy inequality under this modification of the continuity equation, the term $\epsilon (\nabla u_\delta \nabla \rho_\delta)$ is moreover added to the momentum equation. \\
Lastly we have the $\alpha$-level, on which we add the artificial pressure $\alpha \rho^\beta$ is added to the momentum equation. The choice $\beta > \max \lbrace \gamma, 4 \rbrace$ provides us with with a higher regularity of the density, which in turn allows us to pass to the limit in the quantity $\epsilon (\nabla u_\delta \nabla \rho_\delta)$ in the limit passage with respect to $n \rightarrow \infty$, see \cite[Section 7.8.2]{novotnystraskraba}.

\subsection{Solution to the approximate problem} \label{solapproxprob}

Our proof of the existence of a solution to the approximate problem \eqref{approxrho}--\eqref{initialcondapproxi} mainly follows the classical existence theory for the compressible Navier-Stokes equations (c.f.\ for example \cite[Section 7.7]{novotnystraskraba}) with the difference lying only in the consideration of the additional boundary integral in the momentum equation \eqref{weakmomeq}. We start by fixing an arbitrary function $w \in C([0,T];V_n)$. Then by the classical theory for the parabolic Neumann problem (see \cite[Lemma 3.1, Theorem 10.22, Theorem 10.23]{singularlimits}, \cite[Proposition 7.39]{novotnystraskraba}) there exists a unique function $\rho = \rho(w) \in W$ which solves the problem
\begin{align}
&\quad \quad \quad \quad \partial_t \rho + \nabla \cdot \left( \rho w \right) = \epsilon \Delta \rho \quad \text{in } [0,T] \times \Omega, \label{-228} \\
&\rho(0, \cdot) = \rho_0(\cdot)\quad \text{in } \Omega,\ \quad \quad 0 < \underline{\rho} \leq \rho_0(\cdot) \leq \overline{\rho} < \infty, \quad \text{in } \Omega \label{-231}
\end{align}

and which in addition satisfies the estimate
\begin{align}
0 < \underline{\rho} \exp \left(-\left\| w \right\|_{L^1(0,t;V_n)} \right) \leq \rho(t,\cdot) \leq \overline{\rho} \exp \left(\left\| w \right\|_{L^1(0,t;V_n)} \right) < \infty \quad \text{in } \overline{\Omega} \label{-233}
\end{align}

for all $t \in [0,T]$. Further this solution satisfies the estimates
\begin{align}
\left\| \rho(w) \right\|_{C \left([0,T];C^{2,\eta}\left(\overline{\Omega}\right)\right)} + \left\| \rho(w) \right\|_{C^1 \left([0,T];C^{0,\eta}\left(\overline{\Omega}\right)\right)} \leq& c(w), \label{-330} \\
\left\| \rho \left(w^1\right) - \rho \left(w^2\right) \right\|_{C([0,T];L^2(\Omega))} \leq& c\left(w^1,w^2 \right) \left\| w^1 - w^2 \right\|_{C([0,T];W^{1,\infty}(\Omega))} \label{-331}
\end{align}

for all $w,w^1,w^2 \in C([0,T];V_n)$, where the constants $c(w), c(w^1,w^2)>0$ are bounded as long as $w,w^1,w^2$ are bounded in the norm on $C([0,T];V_n)$. Moreover, due to the bound \eqref{-233} of $\rho(w)$ away from 0, it is easy to see from the classical theory of ordinary differential equations that there exists a unique solution $u = u(w) \in C([0,T];V_n)$ to the associated linearized problem
\begin{align}
\int_\Omega \partial_t \left(\rho(w) u \right) \cdot \phi \ dx =& \int _\Omega \left( \rho (w) w \otimes u \right): \nabla \phi- 2\nu \mathbb{D}\left(u\right):\mathbb{D} (\phi) - \lambda (\nabla \cdot u ) (\nabla \cdot \phi) \nonumber \\
& + \left(a\rho(w)^\gamma + \alpha \rho(w)^\beta \right) \nabla \cdot \phi + \rho(w) f\cdot \phi - \epsilon \left( \nabla u \nabla \rho(w) \right) \cdot \phi\ dx \nonumber \\
&- \int_{\partial \Omega} g \operatorname{grad} j_\delta(w) \cdot \phi\ d\Gamma \quad \quad \text{in } [0,T], \label{linmomeq} \\
u(0,\cdot) =& u_0(\cdot) \quad \quad \text{in } \Omega. \label{linincond}
\end{align}

This allows us to consider the desired solution $u_\delta$  to the momentum equation \eqref{weakmomeq} as a fixed point of the operator
\begin{align}
\mathbb{T}: C\left([0,T];V_n\right) \rightarrow C\left([0,T];V_n\right),\quad \mathbb{T}(w) := u(w), \nonumber
\end{align}

mapping $w \in C([0,T];V_n)$ to the corresponding solution to the linearized problem \eqref{linmomeq}. The existence of such fixed point follows from the version \cite[Section 9.2.2, Theorem 4]{evans} of the Schauder fixed point theorem. We show that $\mathbb{T}$ is continuous, compact and fixed points of $s\mathbb{T}$ are bounded in $C([0,T];V_n)$ uniformly with respect to $s \in [0,1]$. To this end we introduce the operator
\begin{align}
\mathcal{M}_{\rho (w)(t)}: V_n \rightarrow V_n^*,\quad \left\langle \mathcal{M}_{\rho (w)(t)}v,\phi \right\rangle_{V_n^* \times V_n} := \int_\Omega \rho (w)(t) v \cdot \phi\ dx \quad \quad \forall \phi, v \in V_n. \nonumber
\end{align}

The bound \eqref{-233} of $\rho(w)$ away from zero implies the existence of an inverse $\mathcal{M}_{\rho (w)(t)}^{-1}$ of $\mathcal{M}_{\rho (w)(t)}$, with the properties
\begin{align}
\left\| \mathcal{M}^{-1}_{\rho (w)(t)} \right\|_{\mathcal{L}(V_n^*,V_n)} &\leq \frac{1}{\inf_{(0,T)\times \Omega}\rho(w)}, \label{18} \\
\left\| \mathcal{M}^{-1}_{\rho (w^1)(t)} - \mathcal{M}^{-1}_{\rho (w^2)(t)} \right\|_{\mathcal{L}(V_n^*,V_n)} &\leq \frac{c(n)}{\left(\inf_{(0,T)\times \Omega} \min \left\lbrace \rho \left(w^1\right), \rho \left(w^2\right) \right\rbrace \right)^2} \left\| \rho\left(w^1\right)(t) - \rho \left(w^2\right)(t) \right\|_{L^1(\Omega)}, \label{17}
\end{align}

as well as
\begin{align}
&\partial_t \left\langle \mathcal{M}_{\rho (w)(t)}v(t),\phi \right\rangle_{V_n^* \times V_n} \nonumber \\
=& \left\langle \mathcal{M}_{\rho (w)(t)}^{-1}\mathcal{M}_{\partial_t \rho (w)(t)}\mathcal{M}_{\rho (w)(t)}^{-1}v(t) + \mathcal{M}_{\rho (w)(t)}^{-1} \partial_t v(t),\phi \right\rangle_{V_n^* \times V_n} \quad \text{in } \mathcal{D}'(0,T), \label{timederid} \\
&\left\| \mathcal{M}_{\rho (w)(t)}^{-1} \mathcal{M}_{\partial_ t\rho (w)(t)}\mathcal{M}_{\rho (w)(t)}^{-1} \right\|_{\mathcal{L}(V_n^*,V_n)} \leq \frac{c(n)}{\left(\inf_{(0,T)\times \Omega}\rho(w)\right)^2} \left\| \partial_t \rho (w)(t) \right\|_{L^1(\Omega)} \label{22}
\end{align}

for any $t \in [0,T]$, any $w,w^1,w^2,v \in C([0,T];V_n)$ and any $\phi \in V_n$, c.f.\ \cite[Section 7.7.1]{novotnystraskraba}. Denoting
\begin{align}
\left\langle \mathcal{N}\left(w,\rho, u \right), \phi \right\rangle_{V_n^*\times V_n} :=& \int _\Omega \left( \rho (w) w \otimes u \right): \mathbb{D} (\phi) + a\rho ^\gamma (w) \nabla \cdot \phi + \alpha \rho ^\beta (w) \nabla \cdot \phi \nonumber \\
&- 2\nu \mathbb{D}\left(u \right):\mathbb{D} (\phi) - \lambda (\nabla \cdot u)(\nabla \cdot \phi) + \rho (w) f\cdot \phi - \epsilon \left( \nabla u \nabla \rho (w) \right) \cdot \phi\ dx \nonumber \\
&- \int_{\partial \Omega} g \operatorname{grad} j_\delta (w) \cdot \phi\ d\Gamma , \nonumber \\
\left\langle \left( \rho_0 u_0 \right)^*, \phi \right\rangle_{V_n^*\times V_n} :=& \int _\Omega \rho_0 u_0 \cdot \phi \ dx, \nonumber
\end{align}

for any $\phi \in V_n$, the solution $u = \mathbb{T}(w)$ to the linearized problem \eqref{linmomeq}, \eqref{linincond} can be expressed as
\begin{align}
u(t) = \mathcal{M}_{\rho (w)(t)}^{-1} \left[ \left( \rho_0 u_0 \right)^* + \int_{0}^t \mathcal{N}\left(w,\rho (w), u \right)\ d\tau \right]. \label{15}
\end{align}

Combining this identity with the estimates \eqref{-233}, \eqref{-330}, \eqref{18}, \eqref{17} and the local Lipschitz-continuity of $\operatorname{grad} j_\delta$ we deduce that the operator $\mathbb{T}$ is continuous. Further, the combination of the identity \eqref{15} with the identity \eqref{timederid} and the estimates \eqref{-233}, \eqref{-331}, \eqref{18} and \eqref{22} leads to the estimate
\begin{align}
\left\| \partial_t u \right\|_{L^2(0,T;V_n)}^2 \leq c(n,w) \label{timederest}
\end{align}

with a constant $c(n,w)>0$ which remains bounded as long as $w$ is bounded in the norm of $C([0,T];V_n)$. From this estimate we infer that the operator $\mathbb{T}$ is also compact. Finally we consider an arbitrary number $s \in [0,1]$ and an arbitrary fixed point $u \in C([0,T];V_n)$ of the operator $s\mathbb{T}$. We test the corresponding linearized momentum equation \eqref{linmomeq} by $u$ and subtract from it the corresponding continuity equation \eqref{-228}, tested by $\frac{1}{2}|u|^2$. This yields the energy equality
\begin{align}
&\frac{d}{dt} \int_\Omega \frac{1}{2} \rho (u)(t)|u(t)|^2 + s \frac{a\rho ^\gamma(u)(t)}{\gamma -1} + s \frac{\alpha \rho ^\beta(u)(t)}{\beta - 1}\ dx + \int_\Omega 2\nu |\mathbb{D} (u)(t)|^2 + \lambda \left| \nabla \cdot u \right|^2\ dx \nonumber \\
&+ s\epsilon \gamma \int_\Omega \rho ^{\gamma-2}(u) \left| \nabla \rho (u) \right|^2\ dx + s\epsilon \beta \int_\Omega \rho ^{\beta-2}(u) \left| \nabla \rho (u) \right|^2\ dx + \int_{\partial \Omega} sg \operatorname{grad} j_\delta(u) \cdot u\ d\Gamma \nonumber \\
=& s\int_\Omega \rho (u)(t)f(t) \cdot u(t)\ dx \label{25}
\end{align}

for all $t \in [0,T]$. According to the property \eqref{nonnegj} of $j_\delta$ it holds that $\operatorname{grad} j_\delta(u) \cdot u \geq 0$ and, by assumption, $g$ is nonnegative. Hence, from the Gronwall Lemma, we deduce that all fixed points $u$ of $s\mathbb{T}$ are bounded in the norm of $C([0,T];V_n)$, independently of $s$. This and the continuity as well as the compactness of $\mathbb{T}$ provides the conditions for the fixed point theorem \cite[Section 9.2.2, Theorem 4]{evans}, which proves the existence of a fixed point $u \in C([0,T];V_n)$ of $\mathbb{T}$. Setting $(\rho_\delta,u_\delta)=(\rho(u),u)$, the pair $(\rho_\delta,u_\delta)$ constitutes the desired solution to our approximate problem \eqref{approxrho}--\eqref{initialcondapproxi}. Integrating the energy inequality \eqref{25}, which $\rho_\delta, u_\delta$ satisfy for $s=1$ we have shown the following proposition:
\begin{proposition} 
\label{System on delta level}
Let the conditions of Theorem \ref{mainresult} be satisfied, let $n \in \mathbb{N}$, $\delta , \epsilon, \alpha > 0$ and let $\beta > \max\lbrace 4, \gamma \rbrace$. Moreover, let $u_0 \in V_n$ be defined by \eqref{indata} and assume $\rho_0$, $q$, defined by \eqref{initialdata}, to satisfy the additional regularity conditions \eqref{-342}. Then there exists a solution $(\rho_\delta, u_\delta) \in W \times C([0,T];V_n)$ to the approximate problem \eqref{approxrho}--\eqref{initialcondapproxi} which in addition satisfies the energy equality
\begin{align}
&\int_\Omega \frac{1}{2} \rho_\delta (\tau)|u_\delta(\tau)|^2 + \frac{a\rho_\delta ^\gamma(\tau)}{\gamma -1} + \frac{\alpha \rho_\delta ^\beta(\tau)}{\beta - 1}\ dx + \int_0^\tau \int_\Omega 2\nu |\mathbb{D} \left(u_\delta\right)(t)|^2 + \lambda \left| \nabla \cdot u_\delta \right|^2 + \epsilon \gamma \rho_\delta ^{\gamma-2} \left| \nabla \rho_\delta \right|^2 \nonumber \\
&+ \epsilon \beta \rho_\delta ^{\beta-2} \left| \nabla \rho_\delta \right|^2\ dxdt + \int_0^\tau \int_{\partial \Omega} g \operatorname{grad} j_\delta \left(u_\delta\right) \cdot u_\delta\ d\Gamma dt \nonumber \\
=& \int_0^\tau \int_\Omega \rho_\delta f \cdot u_\delta\ dxdt +  \int_\Omega \frac{1}{2} \rho_0|u_0|^2 + \frac{a\rho_0 ^\gamma}{\gamma -1} + \frac{\alpha \rho_0 ^\beta}{\beta - 1}\ dx\label{approxen}
\end{align}

for all $\tau \in [0,T]$.
\end{proposition}

\subsection{Limit passage with respect to \texorpdfstring{$\delta \rightarrow 0$}{}} \label{rothelimit} \par

Our next goal is to pass to the limit in the regularization of the function $j$, i.e. the approximation parameter $\delta$ tend to zero. From the energy inequality \eqref{approxen}, the equivalence of norms on the finite dimensional space $V_n$ and the estimates \eqref{-233}, \eqref{-330} for the solution $\rho_\delta$ to the Neumann problem \eqref{-228}, \eqref{-231} with $w = u_\delta$, we infer the uniform bounds
\begin{align}
\left\| \rho_\delta \right\|_{C \left([0,T];C^{2,\eta}\left(\overline{\Omega}\right)\right)} + \left\| \rho_\delta \right\|_{C^1 \left([0,T];C^{0,\eta}\left(\overline{\Omega}\right)\right)} + \left\| \frac{1}{\rho_\delta} \right\|_{C([0,T]\times \overline{\Omega})} + \left\| u_{\delta} \right\|_{C([0,T];V_n)} \leq c \nonumber
\end{align}

with a constant $c>0$ independent of $\delta$. In particular, the bound for $u_\delta$ implies that the bound \eqref{timederest} for $\partial_t u_\delta$ still holds true,
\begin{align}
\left\| \partial_t u \right\|_{L^2(0,T;V_n)}^2 \leq c \label{timederestun}
\end{align}

with a constant $c>0$ independent of $\delta$. Consequently, making use of the Aubin-Lions Lemma, we may extract subsequences and find functions
\begin{align}
0 \leq \rho \in& \bigg\lbrace \psi \in C\left([0,T];H^{1,2}(\Omega)\right) \bigcap C\left([0,T];L^p(\Omega) \right) \bigcap L^2 \left(0,T; H^{2,2}(\Omega) \right): \nonumber \\
&\ \ \partial_t \psi \in L^2 \left((0,T)\times \Omega \right),\ \left. \nabla \psi \cdot \operatorname{n} \right|_{\partial \Omega} = 0 \bigg\rbrace \quad \quad \forall 1\leq p < \infty, \label{-354} \\
u \in& \left\lbrace \phi \in C\left([0,T];V_n \right):\ \partial_t \phi \in L^2\left(0,T;V_n \right) \right\rbrace, \nonumber
\end{align}

such that
\begin{align}
\rho_{\delta} \rightarrow& \rho \quad \text{in } C\left([0,T];H^{1,2}(\Omega) \right) \ \ \text{and} \ \ C\left([0,T];L^p(\Omega) \right), \quad \quad u_{\delta} \rightarrow u \ \ \quad \text{in } C\left([0,T];V_n \right), \label{deltacon1} \\
\rho_{\delta} \rightharpoonup \rho \quad &\text{in } L^2\left(0,T;H^{2,2}(\Omega) \right),\quad \quad \partial_t \rho_{\delta} \rightharpoonup \partial_t \rho \quad \text{in } L^2\left((0,T)\times \Omega \right),\quad \quad \partial_t u_{\delta} \rightharpoonup \partial_t u \quad \text{in } L^2\left(0,T;V_n \right). \label{deltacon2}
\end{align}

Clearly, these convergences are sufficient to pass to the limit in the continuity equation \eqref{weakconteq} and infer that the limit functions $\rho$, $u$ satisfy
\begin{align}
\partial_t \rho + \nabla \cdot \left( \rho u \right) =& \epsilon \Delta \rho \quad \text{a.e. in } (0,T) \times \Omega. \label{weakconteqn}
\end{align}

As in the weak formulation \eqref{momentumequation} we want to combine the momentum equation and the energy inequality into one single relation. To this end we integrate the momentum equation \eqref{weakmomeq} over $[0,\tau]$ for some arbitrary $\tau \in [0,T]$ and subtract from it the energy equality \eqref{approxen}. Further, we exploit the convexity and the $C^1$-regularity of $j_\delta$ to estimate
\begin{align}
\operatorname{grad} j_\delta \left(u_\delta \right) \cdot \left(\phi - u_\delta \right) \leq j_\delta (\phi) - j_\delta \left(u_\delta \right), \nonumber
\end{align}

which allows us to bring the boundary integrals into the same form as in the weak formulation \eqref{momentumequation}. Altogether we obtain the inequality
\begin{align}
&\int_\Omega \frac{1}{2} \rho_0|u_0|^2 + a \frac{\rho_0^\gamma}{\gamma -1} + \frac{\alpha \rho_0 ^\beta}{\beta - 1}\ dx - \int_\Omega \frac{1}{2} \rho_\delta(\tau)|u_\delta(\tau)|^2 + a \frac{\rho_\delta ^\gamma(\tau)}{\gamma -1} + \frac{\alpha \rho_\delta ^\beta(\tau)}{\beta - 1}\ dx \nonumber \\
&+\int_0^\tau \int_\Omega - \rho_\delta u_{\delta} \cdot \partial_t \phi - \left( \rho_\delta u_\delta \otimes u_{\delta} \right): \nabla \phi +  2\nu \mathbb{D}\left(u_{\delta}\right):\mathbb{D} (\phi - u_\delta) + \lambda (\nabla \cdot u_\delta ) (\nabla \cdot (\phi - u_\delta )) \nonumber \\
& - a\rho_\delta ^\gamma \nabla \cdot \phi - \alpha \rho_\delta ^\beta \nabla \cdot \phi 
 - \epsilon a \rho_\delta^{\gamma - 2} \left| \nabla \rho_\delta \right|^2 - \epsilon \beta \rho_\delta^{\beta - 2} \left| \nabla \rho_\delta \right|^2 + \epsilon \left( \nabla u_{\delta} \nabla \rho_\delta \right) \cdot \phi - \rho_\delta f\cdot (\phi - u_\delta)\ dxdt \nonumber \\
&+ \int_0^\tau \int_{\partial \Omega} g j_\delta (\phi) - gj_\delta (u_\delta ) \ d\Gamma dt \geq 0 \quad \quad \forall \phi \in C_c^1 \left((0,\tau);V_n \right),\ \tau \in [0,T]. \label{combinedeq}
\end{align}

Due to the uniform convergences \eqref{convergencej} of $j_\delta$ and \eqref{deltacon1} of $u_\delta$ we can pass to the limit in the boundary integral,
\begin{align}
\int_0^\tau \int_{\partial \Omega} g j_\delta (\phi) - gj_\delta (u_\delta ) \ d\Gamma dt \rightarrow \int_0^\tau \int_{\partial \Omega} g \left| \phi \right| - g \left| u \right| \ d\Gamma dt. \nonumber
\end{align}

The strong convergences \eqref{deltacon1} also allow us to pass to the limit in the remaining terms of the inequality \eqref{combinedeq}. Hence, dropping the nonpositive quantity $-\epsilon a \rho_\delta^{\gamma -2}|\nabla \rho_\delta|^2$ from the left-hand side of this inequality we conclude that the limit functions $\rho$, $u$ satisfy
\begin{align}
&\int_\Omega \frac{1}{2} \rho_0|u_0|^2 + a \frac{\rho_0^\gamma}{\gamma -1} + \frac{\alpha \rho_0 ^\beta}{\beta - 1}\ dx - \int_\Omega \frac{1}{2} \rho(\tau)|u(\tau)|^2 + a \frac{\rho ^\gamma(\tau)}{\gamma -1} + \frac{\alpha \rho ^\beta(\tau)}{\beta - 1}\ dx \nonumber \\
&+\int_0^\tau \int_\Omega - \rho u \cdot \partial_t \phi - \left( \rho u \otimes u \right): \nabla \phi +  2\nu \mathbb{D}\left(u\right):\mathbb{D} (\phi - u) + \lambda (\nabla \cdot u ) (\nabla \cdot (\phi - u)) \nonumber \\
& - a\rho^\gamma \nabla \cdot \phi - \alpha \rho^\beta \nabla \cdot \phi 
 - \epsilon \beta \rho^{\beta - 2} \left| \nabla \rho \right|^2 + \epsilon \left( \nabla u \nabla \rho \right) \cdot \phi - \rho f\cdot (\phi - u)\ dxdt \nonumber \\
&+ \int_0^\tau \int_{\partial \Omega} g \left| \phi \right| - g \left| u \right| \ d\Gamma dt \geq 0 \quad \quad \forall \phi \in C_c^1 \left((0,\tau);V_n \right),\ \tau \in [0,T]. \label{combinedeqn}
\end{align}

\begin{remark} \label{remarkaltmomeq}
As a technical tool we will  need in  \eqref{cweakmom}, \eqref{improvedpressesteps} and  \eqref{effviscflux} a limit version of the momentum equation \eqref{weakmomeq} itself. In this limit equation it will be sufficient to restrict ourselves to test functions vanishing on $\Gamma$. Using such test functions in the momentum equation \eqref{weakmomeq}, we see that the boundary integral vanishes and we can pass to the limit to obtain the identity 
\begin{align}\label{alternativmomeq}
-\int_0^T\int_\Omega \rho u \cdot \partial_t \phi \ dxdt =& \int_0^T\int _\Omega \left( \rho u \otimes u \right): \nabla \phi - 2\nu \mathbb{D}\left(u\right):\mathbb{D} (\phi) - \lambda (\nabla \cdot u ) (\nabla \cdot \phi) + a\rho ^\gamma  \nabla \cdot \phi\nonumber \\
& + \alpha \rho ^\beta \nabla \cdot \phi + \rho f\cdot \phi - \epsilon \left( \nabla u \nabla \rho \right) \cdot \phi\ dxdt
\end{align}

for all $\phi \in C_c^1((0,T);V_n)$ such that $\phi|_{\Gamma} = 0$.
\end{remark}

\subsection{Limit passage with respect to \texorpdfstring{$n \rightarrow \infty$}{}} \par

In this section we pass to the limit in the Galerkin approximation,\ i.e. we let $n$ tend to infinity. Choosing $\phi = 0$ in the momentum and energy inequality \eqref{combinedeqn} we obtain a classical energy inequality. In combination with the classical regularity for the regularized continuity equation \eqref{weakconteqn} (see for example \cite[Lemma 7.37, Lemma 7.38, Section 7.8.2]{novotnystraskraba}) this yields the uniform bounds
\begin{align}
\left\| \rho_n |u_n|^2 \right\|_{L^\infty (0,T;L^1(\Omega))} + \left\| \rho_n \right\|_{L^\infty (0,T;L^\beta (\Omega))} + \left\| u_n \right\|_{L^2(0,T;H^1(\Omega))} \leq& c, \label{nbounds1} \\
\epsilon^\frac{1}{2} \left\| \rho_n^\frac{\beta}{2} \right\|_{L^2(0,T;H^1(\Omega))} + \epsilon \left\| \nabla \rho_n \right\|_{L^{r}((0,T)\times \Omega)} + \epsilon \left\| \partial_t \rho_n \right\|_{L^{\tilde{r}}((0,T)\times \Omega)} + \epsilon^2 \left\| \Delta \rho_n \right\|_{L^{\tilde{r}}((0,T)\times \Omega)} \leq& c \label{nbounds2}
\end{align}

for a constant $c>0$ independent of $n \in \mathbb{N}$, where we can choose
\begin{align}
r := \frac{10\beta - 6}{3\beta + 3} > 2,\quad \tilde{r} := \frac{5\beta - 3}{4\beta} > 1, \nonumber
\end{align}

provided that $\beta > 4$. Interpolations between these bounds lead to the uniform bounds
\begin{align}
\left\| \rho_n u_n \right\|_{L^\infty (0,T;L^\frac{2\beta}{\beta + 1}(\Omega))} + \left\| \rho_n u_n \otimes u_n \right\|_{L^\frac{6\beta}{4\beta + 3}((0,T)\times \Omega)} + \epsilon^{\frac{3}{5\beta}} \left\| \rho_n \right\|_{L^{\frac{5}{3}\beta}((0,T)\times \Omega)} \leq c \label{nbounds3}
\end{align}

for another constant $c>0$ independent of $n$. Combining the bounds \eqref{nbounds1}--\eqref{nbounds3} with the Aubin-Lions Lemma we may extract a subsequence and conclude the existence of functions $u \in L^2(0,T;H_{\text{n}}^1(\Omega))$, $\overline{\rho u \otimes u} \in L^{\frac{6\beta}{4\beta +3}}((0,T)\times \Omega)$ and
\begin{align}
0 \leq \rho \in& L^\infty\left(0,T;L^\beta(\Omega)\right) \bigcap L^{2} \left(0,T; H^1(\Omega) \right) \bigcap L^{\tilde{r}}\left(0,T; W^{2,\tilde{r}}(\Omega) \right) \nonumber
\end{align}

with the properties
\begin{align}
\partial_t \rho \in L^{\tilde{r}} \left((0,T) \times \Omega \right),\quad \quad \left. \nabla \rho \cdot \operatorname{n} \right|_{\partial \Omega} = 0 \label{limitproperties}
\end{align}

such that
\begin{align}
\rho_n& \rightarrow \rho \quad \text{in } L^\beta \left(Q \right)\ \ \text{and} \ \ L^2\left(0,T;H^{1,2}(\Omega) \right),\quad \quad \rho_n \rightharpoonup \rho \quad \text{in } L^{\tilde{r}}\left(0,T;W^{2,\tilde{r}}(\Omega) \right), \label{nconv1} \\
&\quad \ u_n \rightharpoonup u \quad \text{in } L^2 (0,T;H^{1,2} (\Omega)),\quad \quad \partial_t \rho_n \rightharpoonup \partial_t \rho \quad \text{in } L^{\tilde{r}}\left((0,T)\times \Omega \right), \label{nconv2} \\
\rho_n u_n& \buildrel\ast\over\rightharpoonup \rho u \quad \text{in } L^\infty \left(0,T;L^\frac{2\beta}{\beta + 1}(\Omega)\right),\quad \quad \rho_n u_n \otimes u_n \rightharpoonup \overline{\rho u \otimes u} \quad \text{in } L^\frac{6\beta}{4\beta + 3}\left((0,T)\times \Omega \right). \label{nconv3}
\end{align}

With these convergences, we can directly pass to the limit in the continuity equation \eqref{weakconteqn} and infer that
\begin{align}
\partial_t \rho + \nabla \cdot \left( \rho u \right) =& \epsilon \Delta \rho \quad \text{a.e. in } (0,T) \times \Omega. \label{weakconteqeps}
\end{align}

This equation immediately implies that also the renormalized regularized continuity equation
\begin{align}
\partial_t \zeta \left(\rho \right) + \nabla \cdot \left( \zeta \left(\rho \right)u \right) + \left[ \zeta' \left(\rho \right)\rho - \zeta \left(\rho \right) \right] \nabla \cdot u - \epsilon \Delta \zeta \left(\rho \right) = - \epsilon \zeta'' \left(\rho \right)\left| \nabla \rho \right|^2 \leq 0 \quad \text{a.e. in } (0,T)\times \Omega \label{renconeps}
\end{align}

holds true for all convex functions $\zeta \in C^2([0,+\infty))$. In order to pass to the limit in the momentum and energy inequality \eqref{combinedeqn} we need to identify the limit function $\overline{\rho u \otimes u}$ from the convergence \eqref{nconv3}. To this end we test the alternative momentum equation \eqref{alternativmomeq} by $\psi \phi$ for some arbitrary functions $\psi \in \mathcal{D}(0,T)$ and $\phi \in V_N$, $N \leq n$, such that $\phi|_\Gamma = 0$. Under exploitation of the uniform bounds \eqref{nbounds1}--\eqref{nbounds3} this leads us to the dual estimate
\begin{align}
&\left\| \partial_t \int_\Omega \rho_n u_n \cdot \phi \ dx \right\|_{L^{\min \left\{\frac{6}{5}, \frac{2r}{2+r} \right\}}(0,T)} \leq c \nonumber
\end{align}

for a constant $c>0$ depending on $N$ but not on $n$. This allows us to infer from the Arzel\`{a} - Ascoli theorem that
\begin{align}
\int_\Omega \rho_n(\cdot, x) u_n (\cdot, x) \cdot \phi \ dx \rightarrow \int_\Omega \rho (\cdot, x) u(\cdot, x) \cdot \phi \ dx \quad \text{in } C\left([0,T] \right) \label{finiteconvergence}
\end{align}

for any fixed $\phi \in V_N$, $N \in \mathbb{N}$. Since the Galerkin spaces $V_N$ have been choosen such that the functions $\phi \in \bigcup_{N=1}^\infty V_N$ with $\phi|_\Gamma = 0$ are dense in $L^\frac{2\beta}{\beta +1}(\Omega)$ (c.f.\ \eqref{galerkindensity2}) and due to the continuity of the functions $\rho_n u_n$ with respect to the time variable (c.f.\ \eqref{deltacon1}) the convergence \eqref{finiteconvergence} suffices to infer that
\begin{align}
\rho_n u_n \rightarrow \rho u \quad \text{in } C_{\text{weak}}\left([0,T];L^\frac{5}{4}(\Omega) \right) \ \ \text{and thus in}\ \ L^2\left(0,T;H^{-1}(\Omega) \right), \label{cweakmom}
\end{align}

which is sufficient to identify, as desired,
\begin{align}
\overline{\rho u \otimes u} = \rho u \otimes u \quad \text{a.e. in } (0,T) \times \Omega. \label{convidentity}
\end{align}

For the limit passage in the boundary integrals we note that by the weak convergence of $u_n$ in $L^2(0,T;H^1(\Omega))$ and the trace theorem, $u_n$ also converges weakly in $L^2((0,T)\times \partial \Omega)$. Hence the nonnegativity of $g \in L^2((0,T)\times \partial \Omega)$ and the weak lower semicontinuity of the $L^1((0,T)\times \partial \Omega)$-norm imply that
\begin{align}
\int_0^\tau \int_{\partial \Omega} g \left|u \right| \ d\Gamma dt \leq \liminf_{n \rightarrow \infty} \int_0^\tau \int_{\partial \Omega} g \left|u _n\right| \ d\Gamma dt \quad \quad \forall \tau \in [0,T]. \nonumber
\end{align}

This, in combination with the convergences \eqref{nconv1}--\eqref{nconv3}, the identification \eqref{convidentity} of the weak limit of the convective term and the weak lower semicontinuity of norms, gives us all the ingredients required for passing to the limit in both the momentum and energy inequality \eqref{combinedeq} and the alternative momentum equation \eqref{alternativmomeq}. Due to the density of the Galerkin functions in $W_{\text{n}}^{1,p}(\Omega)$ and $W_0^{1,p}(\Omega)$, $1 \leq p < \infty$ (c.f.\ \eqref{galerkindensity}, \eqref{galerkindensity2}) we infer that
\begin{align}
&\int_\Omega \frac{1}{2} \rho_0|u_0|^2 + a \frac{\rho_0^\gamma}{\gamma -1} + \frac{\alpha \rho_0 ^\beta}{\beta - 1}\ dx - \int_\Omega \frac{1}{2} \rho(\tau)|u(\tau)|^2 + a \frac{\rho ^\gamma(\tau)}{\gamma -1} + \frac{\alpha \rho ^\beta(\tau)}{\beta - 1}\ dx \nonumber \\
&+\int_0^\tau \int_\Omega - \rho u \cdot \partial_t \phi - \left( \rho u \otimes u \right): \nabla \phi +  2\nu \mathbb{D}\left(u\right):\mathbb{D} (\phi - u) + \lambda (\nabla \cdot u ) (\nabla \cdot (\phi - u)) \nonumber \\
& - a\rho^\gamma \nabla \cdot \phi - \alpha \rho^\beta \nabla \cdot \phi 
 - \epsilon \beta \rho^{\beta - 2} \left| \nabla \rho \right|^2 + \epsilon \left( \nabla u \nabla \rho \right) \cdot \phi - \rho f\cdot (\phi - u)\ dxdt \nonumber \\
&+ \int_0^\tau \int_{\partial \Omega} g \left| \phi \right| - g \left| u \right| \ d\Gamma dt \geq 0 \label{combinedeqeps}
\end{align}

holds true for almost all $\tau \in [0,T]$ and all $\phi \in \mathcal{D}((0,\tau)\times \overline{\Omega})$ with $\phi \cdot \text{n} |_{\partial \Omega} = 0$ and
\begin{align}\label{alternativmomeqeps}
-\int_0^T\int_\Omega \rho u \cdot \partial_t \phi \ dxdt =& \int_0^T\int _\Omega \left( \rho u \otimes u \right): \nabla \phi - 2\nu \mathbb{D}\left(u\right):\mathbb{D} (\phi) - \lambda (\nabla \cdot u ) (\nabla \cdot \phi) + a\rho ^\gamma  \nabla \cdot \phi\nonumber \\
& + \alpha \rho ^\beta \nabla \cdot \phi + \rho f\cdot \phi - \epsilon \left( \nabla u \nabla \rho \right) \cdot \phi\ dxdt,
\end{align}

for all $\phi \in \mathcal{D}((0,T)\times \Omega)$.

\subsection{Limit passage with respect to \texorpdfstring{$\epsilon \rightarrow 0$}{}} \par \label{epslimit}

In this section we consider the limit passage with respect to $\epsilon \rightarrow 0$ in order to get rid of the artificial regularization terms in the system. By setting $\phi = 0$ in the momentum and energy inequality \eqref{combinedeqeps} and a subsequent interpolation we infer, exactly as the corresponding bounds \eqref{nbounds1} and \eqref{nbounds3} in the previous limit passage, the uniform bounds
\begin{align}
\left\| \rho_\epsilon |u_\epsilon|^2 \right\|_{L^\infty (0,T;L^1(\Omega))} + \left\| \rho_\epsilon \right\|_{L^\infty (0,T;L^\beta (\Omega))} + \left\| u_\epsilon \right\|_{L^2(0,T;H^1(\Omega))} \leq& c, \label{epsbounds1} \\
\left\| \rho_\epsilon u_\epsilon \right\|_{L^\infty (0,T;L^\frac{2\beta}{\beta + 1}(\Omega))} + \left\| \rho_\epsilon u_\epsilon \otimes u_\epsilon \right\|_{L^\frac{6\beta}{4\beta + 3}((0,T)\times \Omega)} \leq& c \label{epsbounds2}
\end{align}

for a constant $c>0$ independent of $\epsilon$. These bounds allow us to extract a subsequence and conclude the existence of functions
\begin{align}
0 \leq \rho \in L^\infty \left(0,T;L^\beta (\Omega) \right),\quad \quad u\in L^2\left(0,T;H_{\text{n}}^1(\Omega)\right) \label{limitpropertiesdelta}
\end{align}

such that
\begin{align}
\rho_\epsilon \buildrel\ast\over\rightharpoonup \rho \ \ \ \text{in } L^\infty \left(0,T;L^\beta(\Omega)\right), \quad \quad u_\epsilon \rightharpoonup u \ \ \ \text{in } L^2 \left(0,T;H^{1,2}(\Omega)\right). \label{epsconv1}
\end{align}

Under exploitation of the continuity equation \eqref{weakconteqeps} the first one of these convergences further leads to
\begin{align}
\rho_\epsilon \rightarrow \rho \quad \text{in } C_{\text{weak}} \left([0,T];L^\beta (\Omega) \right) \quad \quad \text{and hence} \quad \quad \rho_\epsilon u_\epsilon \rightharpoonup \rho u \quad \text{in } L^\infty \left(0,T;L^\frac{2\beta}{\beta +1}(\Omega) \right). \label{epsconv2}
\end{align}

Moreover, we may test the continuity equation \eqref{weakconteqeps} by $\rho_\epsilon$ to infer that
\begin{align}
\epsilon^\frac{1}{2} \left\| \nabla \rho_\epsilon \right\|_{L^2((0,T)\times \Omega)} \leq c \quad \quad \text{and thus} \quad \quad \epsilon \nabla \rho _\epsilon \rightarrow 0 \quad \text{in } L^2\left((0,T)\times \Omega \right). \label{epsconv3}
\end{align}

The convergences \eqref{epsconv1}--\eqref{epsconv3} allow us to pass to the limit in the continuity equation \eqref{weakconteq} and infer that $\rho$ and $u$ satisfy the continuity equation \eqref{-214} in $\mathcal{D}'((0,T)\times \Omega)$. Due to the Lipschitz regularity of $\partial \Omega$ $u$ can be extended continuously to a function $u \in L^2(0,T;H^1(\mathbb{R}^3))$, see \cite[Section 5.4, Theorem 1]{evans}. Extending also $\rho$ by $0$ outside of $\Omega$ we infer that the continuity equation in fact holds true in $\mathcal{D}'((0,T)\times \mathbb{R}^3)$. From the regularization method by DiPerna and Lions (see \cite[Theorem 6.9]{novotnystraskraba}) it follows that $\rho$ and $u$ also satisfy the renormalized continuity equation \eqref{-215}, \eqref{-216}. Similar to the convergence of the convective term in the Galerkin limit (c.f.\ \eqref{nconv3}, \eqref{convidentity}), we may deduce from the alternative momentum equation \eqref{alternativmomeqeps} that
\begin{align}
\rho_\epsilon u_\epsilon \otimes u_\epsilon \rightharpoonup \rho u \otimes u \quad \text{in } L^\frac{6\beta}{4\beta + 3}\left((0,T)\times \Omega \right). \label{epsconv4}
\end{align}

In order to pass to the limit in the pressure terms we need to find a uniform bound for the densitiy in $L^q((0,T)\times \Omega)$ for some $q>\beta$. Thanks to the alternative momentum equation \eqref{alternativmomeqeps} such  bound can be derived as in the case of the no-slip boundary condition, c.f.\ \cite[Lemma 3.1]{fnp}. Nameley, since the Bogovskii operator $\mathcal{B}_\Omega$ on $\Omega$ (c.f.\ \cite[Section 3.3.1.2]{novotnystraskraba}) maps zero-mean functions in $L^p(\Omega)$, $1 < p < \infty$, into $W_0^{1,p}(\Omega)$, we can test the alternative momentum equation \eqref{alternativmomeqeps} by functions of the form
\begin{align}
\phi_\epsilon(t,x) := \psi(t) \mathcal{B}_\Omega \left[ \rho_\epsilon (t) - \frac{1}{|\Omega|} \int_\Omega \rho_\epsilon (t,y)\ dy \right](x),\quad 0 \leq \psi \in \mathcal{D}\left(0,T \right). \label{testfunctionspressest}
\end{align}

As $\mathcal{B}_\Omega$ can be understood as an inverse to the divergence operator, this allows us to find a constant $c>0$ independent of $\epsilon$ such that
\begin{align}
\left\| \rho_\epsilon \right\|_{L^{\gamma + 1}((0,T) \times \Omega)} + \left\| \rho_\epsilon \right\|_{L^{\beta + 1}((0,T) \times \Omega)} \leq c. \label{improvedpressesteps}
\end{align}

Thus we find subsequences and functions $\overline{\rho^\gamma} \in L^\frac{\gamma + 1}{\gamma} ((0,T)\times \Omega)$, $\overline{\rho^\beta} \in L^\frac{\beta + 1}{\beta} ((0,T)\times \Omega)$ such that
\begin{align}
\rho_\epsilon^\gamma \rightharpoonup \overline{\rho^\gamma} \quad \text{in } L^\frac{\gamma + 1}{\gamma}\left((0,T)\times \Omega)\right),\quad \quad \rho^\beta_\epsilon \rightharpoonup \overline{\rho^\beta} \quad \text{in } L^\frac{\beta + 1}{\beta}\left((0,T)\times \Omega \right). \label{epsconv5}
\end{align}

Our next goal is to identify the limit functions $\overline{\rho^\gamma}$ and $\overline{\rho^\beta}$, for which we need the effective viscous flux identity
\begin{align}
\lim_{\epsilon \rightarrow 0} \int_{0}^T \int_{\Omega} \Phi (\lambda + 2\nu)  \left[ \rho_\epsilon \nabla \cdot u_\epsilon - \rho \nabla \cdot u \right] \ dxdt = \lim_{\epsilon \rightarrow 0} \int_{0}^T \int_{\Omega} \Phi \left( \left[ a\rho_\epsilon^\gamma + \alpha \rho_\epsilon^\beta \right] \rho_\epsilon - \left[ a\overline{\rho^\gamma} + \alpha \overline{\rho^\beta} \right] \rho \right)\ dxdt \label{effviscflux}
\end{align}

for all $0\leq \Phi \in \mathcal{D}((0,T)\times \Omega)$. This identity can be proved by applying the method from \cite[Lemma 3.2]{fnp} to the alternative momentum equation \eqref{alternativmomeqeps}. We test \eqref{alternativmomeqeps} and a corresponding limit identity, obtained from the convergences \eqref{epsconv1}, \eqref{epsconv2} and \eqref{epsconv5}, by functions of the form
\begin{align}
\phi_\epsilon (t,x) := \Phi(t,x) \left( \nabla \Delta^{-1} \right)\left[ \rho_\epsilon (t,\cdot) \right](t,x),\quad \phi (t,x) := \Phi(t,x) \left( \nabla \Delta^{-1} \right)\left[ \rho(t,\cdot) \right](t,x), \label{testfunctions}
\end{align}

with $0 \leq \Phi \in \mathcal{D}((0,T)\times \Omega)$, respectively. Subtracting the two resulting relations from each other we obtain the effective viscous flux identity exactly as in \cite[Lemma 3.2]{fnp}. Following the procedure in \cite[Section 3.5]{fnp} we consider - after a dominated convergence argument - both the renormalized continuity equation \eqref{renconeps} on the $\epsilon$-level and the renormalized continuity equation \eqref{-215} in the limit with the choice of the (strictly) convex function $\zeta (\xi):= \xi \ln (\xi)$ in. This results in two relations which we subtract from each other to obtain the inequality
\begin{align}
\lim_{\epsilon \rightarrow 0}  \int_\Omega \rho(\tau) \ln \left( \rho (\tau) \right)\ - \rho_\epsilon (\tau) \ln \left( \rho_\epsilon (\tau) \right)\ dx \geq \lim_{\epsilon \rightarrow 0} \int_0^\tau \int_\Omega \rho_\epsilon \nabla \cdot u_\epsilon - \rho \nabla \cdot u \ dxdt \quad \quad \forall \tau \in [0,T]. \label{-157}
\end{align}

From the effective viscous flux identity \eqref{effviscflux} and the monotonicity of the (artificial) pressure function it follows that the right-hand side of this relation is nonnegative.  Further, since the mapping $\xi \rightarrow \xi \ln (\xi)$ is convex, we know that $\rho \ln (\rho) \leq \overline{\rho \ln (\rho)}$, where $\overline{\rho \ln (\rho)}$ denotes a weak limit of $\overline{\rho_\epsilon \ln (\rho_\epsilon)}$ in $L^1((0,T)\times \Omega)$. Combining these two facts, we conclude that
\begin{align}
\rho \ln (\rho) = \overline{\rho \ln (\rho)} \quad \quad \text{a.e. in } (0,T) \times \Omega. \nonumber
\end{align}

By the relations between weakly convergent sequences and (strictly) convex functions (c.f.\ \cite[Theorem 10.20]{singularlimits}), this equation implies pointwise convergence of $\rho_\epsilon$, which in turn implies that, as desired,
\begin{align}
\overline{\rho^\gamma} = \rho^\gamma \quad \text{a.e. in } (0,T)\times \Omega,\quad \quad \overline{\rho^\beta} = \rho^\beta \quad \text{a.e. in } (0,T)\times \Omega. \label{pressureidentificationeps}
\end{align}

Combining the convergences \eqref{epsconv1}--\eqref{epsconv4}, \eqref{epsconv5}, the identification \eqref{pressureidentificationeps} of the limits of the pressure terms and the weak lower semicontinuity of norms, we can now pass to the limit in both the momentum and energy inequality \eqref{combinedeqeps} and the alternative momentum equation \eqref{alternativmomeqeps} and infer that
\begin{align}
&\int_\Omega \frac{1}{2} \rho_0|u_0|^2 + a \frac{\rho_0^\gamma}{\gamma -1} + \frac{\alpha \rho_0 ^\beta}{\beta - 1}\ dx - \int_\Omega \frac{1}{2} \rho(\tau)|u(\tau)|^2 + a \frac{\rho ^\gamma(\tau)}{\gamma -1} + \frac{\alpha \rho ^\beta(\tau)}{\beta - 1}\ dx \nonumber \\
&+\int_0^\tau \int_\Omega - \rho u \cdot \partial_t \phi - \left( \rho u \otimes u \right): \nabla \phi +  2\nu \mathbb{D}\left(u\right):\mathbb{D} (\phi - u) + \lambda (\nabla \cdot u ) (\nabla \cdot (\phi - u)) \nonumber \\
& - a\rho^\gamma \nabla \cdot \phi - \alpha \rho^\beta \nabla \cdot \phi 
 - \rho f\cdot (\phi - u)\ dxdt + \int_0^\tau \int_{\partial \Omega} g \left| \phi \right| - g \left| u \right| \ d\Gamma dt \geq 0 \label{combinedeqalpha}
\end{align}

holds true for almost all $\tau \in [0,T]$ and all $\phi \in \mathcal{D}((0,\tau)\times \overline{\Omega})$ with $\phi \cdot \text{n}|_{\partial \Omega} = 0$ and
\begin{align}\label{alternativmomeqalpha}
-\int_0^T\int_\Omega \rho u \cdot \partial_t \phi \ dxdt =& \int_0^T\int _\Omega \left( \rho u \otimes u \right): \nabla \phi - 2\nu \mathbb{D}\left(u\right):\mathbb{D} (\phi) - \lambda (\nabla \cdot u ) (\nabla \cdot \phi) + a\rho^\gamma  \nabla \cdot \phi\nonumber \\
& + \alpha \rho^\beta \nabla \cdot \phi + \rho f\cdot \phi dxdt,
\end{align}

holds true for all $\phi \in \mathcal{D}((0,T)\times \Omega)$.

\subsection{Limit passage with respect to \texorpdfstring{$\alpha \rightarrow 0$}{}} \par\label{alphalim}

Finally it remains to get rid of the artificial pressure term in the momentum equation, i.e.\ to let $\alpha$ tend to zero. In addition, we return from the regularized initial data $\rho_{0,\alpha}$, $q_{\alpha}$ in the approximate problem (c.f.\ \eqref{-342}) to the more general initial data $\rho_0$, $q$ from the main result Theorem \ref{mainresult}. More precisely, as in \cite[Section 4]{fnp}, we choose $\rho_{0,\alpha}$, $q_{\alpha}$ satisfying the relations \eqref{-342} for any fixed $\alpha > 0$ such that
\begin{align}
\rho_{0,\alpha} \rightarrow \rho_0 \quad \ \ \ &\text{in } L^\gamma (\Omega),\quad \quad &\alpha \rho_{0,\alpha}^\beta& \rightarrow 0 \quad \quad \quad \quad \text{in } L^1(\Omega), \label{initialconv1} \\
q_{\alpha} \rightarrow q \quad &\text{in } L^1 (\Omega),\quad \quad &\frac{\left|q_{\alpha}\right|^2}{\rho_{0,\alpha}}& \rightarrow \frac{\left|q\right|^2}{\rho_{0}} \quad \quad \ \ \ \text{in } L^1(\Omega) \label{initialconv2}
\end{align}

for $\alpha \rightarrow 0$. As in the previous limit passages we infer, from the choice $\phi = 0$ in the momentum and energy inequality \eqref{combinedeqalpha} and an ensuing interpolation, the uniform bounds
\begin{align}
\left\| \rho_\alpha |u_\alpha|^2 \right\|_{L^\infty (0,T;L^1(\Omega))} + \left\| \rho_\alpha \right\|_{L^\infty (0,T;L^\gamma (\Omega))} + \left\| u_\alpha \right\|_{L^2(0,T;H^1(\Omega))} \leq& c, \label{alphabounds1} \\
\left\| \rho_\alpha u_\alpha \right\|_{L^\infty (0,T;L^\frac{2\gamma}{\gamma + 1}(\Omega))} + \left\| \rho_\alpha u_\alpha \otimes u_\alpha \right\|_{L^\frac{6\gamma}{4\gamma + 3}((0,T)\times \Omega)} \leq& c \label{alphabounds2}
\end{align}

for a constant $c>0$ independent of $\alpha$. This allows us to find a subsequence as well as functions
\begin{align}
0 \leq \rho \in L^\infty \left(0,T;L^\gamma (\Omega) \right),\quad \quad u\in L^2\left(0,T;H_{\text{n}}^1(\Omega)\right) \label{limitpropertiesfinal}
\end{align}

such that
\begin{align}
\rho_\alpha \buildrel\ast\over\rightharpoonup \rho \quad \text{in } L^\infty \left(0,T;L^\gamma(\Omega)\right), \quad \quad u_\alpha \rightharpoonup u \quad \text{in } L^2 \left(0,T;H^{1,2}(\Omega)\right)\label{alphaconv1}
\end{align}

as well as, under exploitation of the continuity equation \eqref{-214} and the alternative momentum equation \eqref{alternativmomeqalpha},
\begin{align}
\rho_\alpha \rightarrow \rho \quad \text{in } C_{\text{weak}} \left([0,T];L^\gamma (\Omega) \right),\quad \quad \rho_\alpha u_\alpha \rightarrow \rho u \quad \text{in } C_{\text{weak}} \left([0,T];L^\frac{2\gamma}{\gamma + 1} (\Omega) \right) \label{alphaconv2}
\end{align}

and consequently
\begin{align}
\rho_\alpha u_\alpha \otimes u_\alpha \rightharpoonup \rho u \otimes u \quad \text{in } L^\frac{6\gamma}{4\gamma + 3}((0,T)\times \Omega). \label{alphaconv3}
\end{align}

Due to the convergences \eqref{alphaconv2} and the continuity equation \eqref{-214} on the $\alpha$-level, the limit functions $\rho$ and $u$ satisfy the same continuity equation \eqref{-214} in $\mathcal{D}'((0,T)\times \Omega)$. For the limit passage in the pressure terms the derivation of the improved uniform bounds \eqref{improvedpressesteps} of the density on the $\epsilon$-level needs to be modified. More specifically, the derivation of these bounds relies on the fact that on the $\epsilon$-level the density is bounded uniformly in $L^\infty(0,T;L^\beta (\Omega))$, which is not the case in our current situation. As a compensation, the density in the test functions \eqref{testfunctionspressest} needs to be replaced by a suitable smooth approximation of $\rho_\alpha^\theta$ for some sufficiently small value $\theta > 0$. This procedure, which is described in detail in \cite[Section 4.1]{fnp}, leads, by a use of the resulting test functions in the alternative momentum equation \eqref{alternativmomeqalpha}, to the desired improved pressure estimates
\begin{align}
\left\| \rho_\alpha \right\|_{L^{\gamma + \theta}((0,T) \times \Omega)} + \alpha^\frac{1}{\beta + \theta}\left\| \rho_\alpha \right\|_{L^{\beta + \theta}((0,T) \times \Omega)} \leq c \nonumber
\end{align}

with a constant $c>0$ independent of $\alpha$. In particular we may extract a subsequence and find a function $\overline{\rho^\gamma} \in L^\frac{\gamma + 1}{\gamma}((0,T)\times \Omega)$ such that
\begin{align}
\rho_\alpha^\gamma \rightharpoonup \overline{\rho^\gamma} \quad \text{in } L^\frac{\gamma + \theta}{\gamma}\left((0,T)\times \Omega \right),\quad \quad \alpha \rho^\beta_\alpha \rightarrow 0 \quad \text{in } L^\frac{\beta + \theta}{\beta}\left((0,T)\times \Omega\right). \label{alphaconv5}
\end{align}

For the identification of the limit function $\overline{\rho^\gamma}$ we further follow the procedure in \cite[Section 4.3]{fnp} and deduce the following modified version of the effective viscous flux identity \eqref{effviscflux} on the $\epsilon$-level,
\begin{align}
&\lim_{\alpha \rightarrow 0} \int_{0}^T \int_{\Omega} \Phi (\lambda + 2\nu)  \left[ T_k\left(\rho_\alpha \right) \nabla \cdot u_\alpha - \overline{T_k\left(\rho \right)} \nabla \cdot u \right] \ dxdt \nonumber \\
=& \lim_{\alpha \rightarrow 0} \int_{0}^T \int_{\Omega} \Phi \left( a\rho_\alpha^\gamma T_k \left( \rho_\alpha \right) - a\overline{\rho^\gamma}\ \overline{T_k\left(\rho \right)} \right)\ dxdt \label{effviscfluxalpha}
\end{align}

for all $\Phi \in \mathcal{D}((0,T)\times \Omega)$, where $T_k \leq 2k$, $k \in \mathbb{N}$, constitutes a suitable smooth and concave cut-off version of the identity function on $[0,\infty)$ and $\overline{T_k(\rho)}$ denotes a weak limit of $T_k(\rho_\alpha)$ in $L^1((0,T)\times \Omega)$. In the derivation of this identity we again have to make up for the lower integrability of the density as compared to the density on the $\epsilon$-level. This is achieved by replacing the test functions \eqref{testfunctions} used on the $\epsilon$-level by test functions of the form
\begin{align}
\phi_\alpha (t,x) := \Phi(x) \left( \nabla \Delta^{-1} \right) \left[ T_k\left( \rho_\alpha \right) (t,\cdot) \right](x),\quad \quad \phi (t,x) := \Phi(x) \left( \nabla \Delta^{-1} \right) \left[  \overline{ T_k\left( \rho \right)}(t,\cdot) \right](x). \nonumber
\end{align}

where $\Phi \in \mathcal{D}((0,T)\times \Omega)$. In our case, we use these test functions in the alternative momentum equation \eqref{alternativmomeqalpha} on the $\alpha$-level and a corresponding limit identity respectively. Comparing the resulting identities we arrive at the desired effective viscous flux identity \eqref{effviscfluxalpha}, exactly as in the proof of \cite[Lemma 4.2]{fnp}. From this identity, the concavity of $T_k$ and the convexity of $\xi \mapsto \xi^\gamma$ we deduce, exactly as in the proof of \cite[Lemma 4.3]{fnp}, boundedness of the oscillation defect measure,
\begin{align}
\textbf{osc}_{\gamma + 1}\left[\rho_\alpha \rightarrow \rho \right]\left((0,T)\times \mathbb{R}^3 \right) := \sup_{k \geq 1} \left[ \limsup_{\alpha \rightarrow 0} \int_0^T \int_{\mathbb{R}^3} \left| T_k \left(\rho_\alpha \right) - T_k\left( \rho \right) \right|^{\gamma + 1}\ dxdt \right] < \infty. \label{oscbound}
\end{align}

Next we choose $\zeta = T_k$ in the renormalized continuity equation \eqref{-215} on the $\alpha$-level and pass to the limit with respect to $\alpha$. Since $\overline{T_k(\rho)}$ is bounded uniformly, the resulting limit identity can be renormalized under exploitation of the regularization technique by DiPerna and Lions, see \cite[Lemma 6.9]{novotnystraskraba}. Subsequently, using the bound \eqref{oscbound} of the oscillation defect measure, we may let $k$ tend to infinity to infer that $\rho$ and $u$ also satisfy the renormalized continuity equation \eqref{-215}. For the details of this procedure we refer to the proof of \cite[Lemma 4.4]{fnp}. Under exploitation of the dominated convergence theorem we may use the choice
\begin{align}
\zeta (\xi) := \zeta_k(\xi) := \xi \int_1^\xi \frac{T_k(s)}{s^2}\ ds \nonumber 
\end{align}

in both the renormalized continuity equation \eqref{-215} on the $\alpha$-level and in the limit. Comparing the resulting equations to each other and passing to the limit with respect to $\alpha \rightarrow 0$ we infer that
\begin{align}
&\lim_{\alpha \rightarrow  0}\int_{\Omega} \left( \zeta_k(\rho_\alpha) - \zeta_k(\rho) \right)(\tau)\ dx + \lim_{\alpha \rightarrow 0} \int_0^\tau \int_{\Omega} T_k\left(\rho_\alpha \right) \nabla \cdot u_\alpha - \overline{T_k(\rho)} \nabla \cdot u \ dxdt \nonumber \\
=&\int_0^\tau \int_{\Omega} T_k(\rho) \nabla \cdot u - \overline{T_k(\rho)} \nabla \cdot u \ dxdt \label{-201}
\end{align}

for all $\tau \in [0,T]$. Here, the second term on the left-hand side is nonnegative, which follows from the effective viscous flux identity \eqref{effviscfluxalpha}, the fact that both the mappings $\xi \mapsto \xi^\gamma$ and $\xi \mapsto T_k(\xi)$ are nondecreasing and the classical relations between weakly convergent sequences and monotone functions (c.f.\ \cite[Theorem 10.19]{singularlimits}). Moreover, the right-hand side of the equation \eqref{-201} vanishes for $k \rightarrow \infty$ as can be seen from the bound \eqref{oscbound} of the oscillation defect measure. Consequently, letting $k$ tend to infinity also on the left-hand side of this identity we infer that
\begin{align}
\lim_{\alpha \rightarrow 0}  \int_\Omega \rho(\tau) \ln \left( \rho (\tau) \right)\ - \rho_\alpha (\tau) \ln \left( \rho_\alpha (\tau) \right)\ dx \geq 0. \nonumber
\end{align}

Exactly as in the limit passage with respect to $\epsilon \rightarrow 0$ (c.f.\ \eqref{pressureidentificationeps}), this estimate yields pointwise convergence of $\rho_\alpha$ and consequently the identity $\overline{\rho^\gamma} = \rho^\gamma$ almost everywhere in $(0,T)\times \Omega$. Therefore, using the convergences \eqref{alphaconv1}--\eqref{alphaconv3} and \eqref{alphaconv5} as well as the weak lower semicontinuity of norms we may pass to the limit in the momentum and energy inequality \eqref{combinedeqalpha} and infer that $\rho$ and $u$ satisfy the momentum and energy inequality \eqref{momentumequation}. Finally we note that $\rho$, as a solution to the renormalized continuity equation \eqref{-215}, is an element of the space $C([0,T];L^1(\Omega))$, c.f.\ \cite[Proposition 4.3]{feireisldynamics}. Due to this continuity in the time variable and the convergence \eqref{initialconv1} of the initial data it satisfies the initial condition $\rho = \rho_0$ in the classical sense. The initial condition stated for $\rho u$ in \eqref{initialcond} follows from the convergences \eqref{initialconv2} and \eqref{alphaconv2} This concludes the proof of Theorem \ref{mainresult}.

\section*{Acknowledgment}
\v S. Ne\v casov\'a  and J. Ogorzaly have been supported by  Praemium Academiae of \v S. Ne\v {c}asov\'{a}. {J. Ogorzaly (first version of paper) was supported by  by the Czech Science Foundation (GA\v CR) through project 19-04243S.}
Further, the work has been supported by the Czech Science Foundation (GA\v CR) through projects GC22-08633J, (for \v S. Ne\v casov\'a  and J. Scherz).
The Institute of Mathematics, CAS is supported by RVO:67985840.

\end{document}